\documentclass{article}

\date{}
\usepackage{amsmath}
\usepackage{amssymb}
\usepackage{amsfonts}

\title{\Large\bf Boundedness properties in function spaces}

\author{{\bf L'ubica Hol\'a}\\
Mathematical Institute, Slovak Academy of Sciences,\\
\v Stef\'anikova 49, SK-814 73 Bratislava, Slovakia\\
{\sf lubica.hola@mat.savba.sk}\\\\
{\bf Ljubi\v sa D.R. Ko\v cinac}\footnote{Corresponding author}\\
University of Ni\v s, Faculty of Sciences and Mathematics, \\18000 Ni\v s, Serbia\\
{\sf lkocinac@gmail.com} }

\newtheorem{theorem}{{\bf Theorem}}[section]

\newtheorem{remark}[theorem]{{\bf Remark}}
\newtheorem{note}[theorem]{{\bf Note}}

\newtheorem{example}[theorem]{{\bf Example}}

\newcommand{\proof}{{\bf Proof.  }}
\newcommand{\eproof}{{\square}}

\begin{document}

\maketitle

\begin{abstract}
Some boundedness properties of function spaces (considered as
topological groups) are studied.
\end{abstract}



\medskip
\section{Introduction}

Let $X$ be a Tychonoff space and $C(X)$ be the set of continuous
real-valued functions defined on $X$. This set can be equipped with
a variety of (natural) topologies: the topology $\tau_p$ of
pointwise convergence, the compact-open topology $\tau_k$, the
topology $\tau_u$ of uniform convergence, the $m$-topology $\tau_m$
\cite{Hew} (sometimes called also the fine topology, Whitney
topology, Morse topology \cite{Hola-Mccoy}), the graph topology
$\tau_{\gamma}$ \cite{naimpally}, and some other. We write $(C(X),
\tau_p) = C_p(X)$, $(C(X), \tau_k) = C_k(X)$, $(C(X), \tau_u) =
C_u(X)$, $(C(X), \tau_m) = C_m(X)$, and $(C(X),\tau_{\gamma}) =
C_{\gamma}(X)$. In this paper we discuss some boundedness properties
of these topological spaces considering them as Hausdorff
topological groups under pointwise addition. The local base of a
function $f\in C_p(X)$, $f\in C_k(X)$, $f\in C_u(X)$, $f\in C_m(X)$
and $f\in C_{\gamma}(X)$, respectively, consists of all sets of the
form
\[
W(f,F,\varepsilon)= \{g\in C_p(X): |g(x)-f(x)| < \varepsilon \mbox{
for each } x\in F\},
\]
\[
W(f,K,\varepsilon)= \{g\in C_k(X): |g(x)-f(x)| < \varepsilon \mbox{
for each } x\in K\},
\]
\[
W(f,\varepsilon)= \{g\in C_u(X): |g(x)-f(x)| < \varepsilon \mbox{
for each } x\in X\},
\]
\[
W(f,\varepsilon)= \{g\in C_m(X): |g(x)-f(x)| < \varepsilon(x) \mbox{
for each } x\in X\},
\]
\[
W(f,\eta) = \{g \in C_{\gamma}(X):|g(x)-f(x)| < \eta(x) \mbox{ for
each } x\in X\},
\]
where $\varepsilon > 0$, $F$ is a finite, and $K$ a compact set in
$X$, $\varepsilon(x)\in C(X)$ is a strictly positive function, and
$\eta(x)$ is a strictly positive lower semicontinuous real-valued
function on $X$.

Clearly, $\tau_p \le \tau_k\le \tau_u\le \tau_m \le \tau_{\gamma}$.

Because all these spaces are homogeneous, when we study their local
properties it is enough to consider a local base in the constantly
zero function $f_0$ (defined by $f_0(x) = 0$ for each $x\in X$)
which is the identity element in each of these topological groups.

Recall that a topological group $(G,\cdot,\tau)$ is said to be:

\begin{itemize}
\item[{\rm (1)}] \emph{$\aleph_0$-bounded} (called also \emph{$\omega$-narrow} \cite{ArhTk})
if for each neighborhood $U$ of the neutral element $e\in G$ there
is countable set $A\subset G$ such that $G = A\cdot U$ \cite{Gu};

\item[{\rm (2)}] \emph{Menger-bounded} or shortly \emph{{\sf
M}-bounded} if for each sequence $(U_n:n\in\mathbb N)$ of
neighborhoods of $e\in G$ there is a sequence $(A_n:n\in\mathbb N)$
of finite subsets of $G$ such that $G = \bigcup_{n\in\mathbb
N}A_n\cdot U_n$ (\cite{Ko1} and independently \cite{Tk} under the
name \emph{$o$-bounded}; see also \cite{BKS, Ko2});

\item[{\rm (3)}] \emph{Hurewicz bounded} or shortly \emph{{\sf
H}-bounded} if for each sequence $(U_n : n\in\mathbb N)$ of
neighborhoods of $e \in G$ there is a sequence $(A_n : n\in\mathbb
N)$  of finite subsets of $G$ such that each $x \in G$ belongs to
all but finitely many sets $A_n \cdot U_n$ \cite{Ko1, Ko2};

\item[{\rm (4)}] \emph{Rothberger bounded} or shortly \emph{\sf
R}-bounded if for each sequence $(U_n:n\in\mathbb N)$ of
neighborhoods of $e\in G$ there is a sequence $(x_n:n\in\mathbb N)$
of elements of $G$ such that $G = \bigcup_{n\in\mathbb N}x_n\cdot
U_n$ \cite{Ko1, Ko2};

%\item[{\rm (5)}] \emph{Gerlits-Nagy bounded} (shortly {\sf
%GN}-bounded) if for each sequence $(U_n:n\in\mathbb N)$ of
%neighborhoods of $e\in G$ there is a sequence $(x_n:n\in\mathbb N)$
%of elements of $G$ such that each $x\in G$ belongs to all but
%finitely many sets $x_n\cdot U_n$ \cite{Ko1, Ko2}.
\end{itemize}

To each of the above properties one can correspond a game on $G$.
Let us demonstrate this for {\sf M}-boundedness. Two players, ONE
and TWO, play a round for each $n\in \mathbb N$. In the $n$-th round
player ONE chooses a basic neighborhood $U_n$ of $e\in G$, and TWO
responds by choosing a finite set $A_n$ in $G$. TWO wins a play
\[
U_1, A_1; U_2, A_2;  \cdots; U_n, A_n; \cdots
\]
if $G = \bigcup_{n\in\mathbb N}A_n\cdot U_n$; otherwise ONE wins.

\smallskip
A topological group $(G,\cdot,\tau)$ is said to be \emph{strictly
{\sf M}-bounded} if TWO has a winning strategy in the above game
\cite{He, Ko1}.

A similar scheme is applied for definitions of games associated to
other boundedness properties. It is clear what means a
\emph{strictly {\sf H}-bounded} and \emph{strictly {\sf R}-bounded}
topological group.

We have the following diagram for relations among the properties we
mentioned for a topological group.
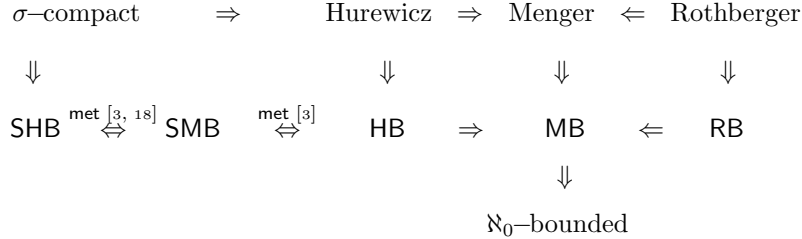
\[
\hskip-0.8cm \sigma\!\!-\!\!\mbox{compact} \hskip0.9cm  \Rightarrow
\hskip0.9cm  \mbox{ Hurewicz } \ \Rightarrow \ \mbox{ Menger } \
\Leftarrow \ \mbox{ Rothberger }
\]
\[
\hskip-1.6cm\Downarrow \hskip4.5cm  \Downarrow \hskip2.1cm
\Downarrow \hskip2cm \Downarrow
\]
\[
\hskip-1.7cm{\sf SHB} \stackrel{{\tiny {\sf met}\ {\cite{Bk, MT}}}}
\Leftrightarrow {\sf SMB} \hskip0.4cm \stackrel{{\tiny {\sf met}\
\cite{Bk}}} \Leftrightarrow \hskip0.6cm {\sf HB} \hskip0.6cm
\Rightarrow \hskip0.7cm  {\sf MB} \hskip0.6cm \Leftarrow \hskip0.5cm
{\sf RB}
\]
\[
\hskip3.3cm\Downarrow
\]
\[
\hskip3.1cm\aleph_0\!\!-\!\!\mbox{bounded}
\]

\medskip
\noindent where {\sf MB}, {\sf SMB}, {\sf HB}, {\sf SHB}, {\sf RB},
and {\sf met} denote {\sf M}-bounded, strictly {\sf M}-bounded, {\sf
H}-bounded, strictly {\sf H}-bounded, {\sf R}-bounded, and
metrizable, respectively.

\medskip
Also, if $\mathcal P \in \{\aleph_0\!\!-\!\!{\rm bounded}, {\sf HB},
{\sf SHB}, {\sf MB}, {\sf SHB},  {\sf RB}, \}$, then for a Tychonoff
space $X$ we have
\[
C_{\gamma}(X) \in \mathcal P \, \Rightarrow \, C_m(X) \in \mathcal P
\, \Rightarrow \, C_u(X) \in \mathcal P \, \Rightarrow  C_k(X) \in
\mathcal P \, \Rightarrow C_p(X) \in \mathcal P.
\]

For some results on {\sf M}-bounded, strictly {\sf M}-bounded, {\sf
H}-bounded and {\sf R}-bounded topological groups we refer the
reader to the papers \cite{Bk, BKS, Ba, HRT, MT, Sch, Tsaban} and
references therein.

\section{Results}

We begin with some basic facts which are known or can be easily
proved.

\smallskip
\noindent {\bf Fact 1.} The class of $\aleph_0$-bounded groups is
productive, hereditary, preserved by continuous epimorphisms
\cite{Gu} and contains the classes of compact, pseudocompact,
totally bounded and separable (thus also ccc) groups (see
\cite[Section 3.4]{ArhTk}). By a result of Guran \cite{Gu}, a
topological group $G$ is $\aleph_0$-bounded if and only if $G$ is
topologically isomorphic to a subgroup of a product of second
countable topological groups.

\smallskip
\noindent {\bf Fact 2.} All the properties from the above diagram
are preserved under continuous epimorphisms.

\smallskip
\noindent {\bf Fact 3.} It is clear that every {\sf M}-bounded group
is $\aleph_0$-bounded. The countable power of $\mathbb R$ is an
$\aleph_0$-bounded group which is not {\sf M}-bounded (see
\cite[Example 2.6]{He}).

\smallskip
\noindent {\bf Fact 4.} The topological group $C_p(X)$ is always
$\aleph_0$-bounded. It follows from the fact that $C_p(X)$ is always
a ccc space \cite{Ar}.

\smallskip
\noindent {\bf Fact 5.} $C_p(X)$ is never {\sf R}-bounded. Suppose
to the contrary that $C_p(X)$ is {\sf R}-bounded. Let $x \in X$.
Then the mapping $\pi: C_p(X) \to C_p(\{x\})$, $f \mapsto
f\restriction \{x\}$, is a continuous epimorphism. Then by Fact 2,
$C_p(\{x\})$ must be {\sf R}-bounded. However, $C_p(\{x\}) \cong
\mathbb R$ and $\mathbb R$ is not {\sf R}-bounded. A contradiction.

\smallskip
\noindent {\bf Fact 6.} If $X$ is a submetrizable space (i.e. $X$
admits a continuous bijection onto a metrizable space), then
$C_k(X)$ is $\aleph_0$-bounded because in this case $C_k(X)$ is a
ccc space \cite[Proposition 7.1.3]{ntantu}.

\smallskip
\noindent {\bf Fact 7.} It is clear that every $\sigma$-compact
group is strictly {\sf H}-bounded and that every subgroup of a
(strictly) {\sf H}-bounded group is (strictly) {\sf H}-bounded. Also
it is shown in \cite{He} that every $\sigma$-compact group is
strictly {\sf M}-bounded and that every subgroup of a (strictly)
{\sf M}-bounded group is (strictly) {\sf M}-bounded.

\begin{remark} \rm (Folklore) Fact 5 follows also from the following claim: a
subgroup $H$ of an {\sf R}-bounded group $G$ is also {\sf
R}-bounded. Let $(U_n:n\in \mathbb N)$ be a sequence of
neighborhoods of the identity element $e \in H$. Take for each $n$ a
neighborhood $V_n$ of $e\in G$ with $V_n\cap H = U_n$. Let for each
$n\in\mathbb N$, $W_n$ be a neighborhood of $e\in G$ such that
$W_n^{-1}\cdot W_n \subset V_n$. Apply to the sequence $(W_n:n\in
\mathbb N)$ the fact that $G$ is {\sf R}-bounded and find a sequence
$(x_n:n\in \mathbb N)$ of elements of $G$ so that $G=
\bigcup_{n\in\mathbb N}x_n\cdot W_n$. Whenever $x_n\cdot W_n \cap H
\neq \emptyset$ pick a point $a_n$ in this intersection; otherwise
take $a_n= e$. The sequence $(a_n:n\in \mathbb N)$ of elements of
$H$ witnesses for $(U_n:n\in \mathbb N)$ that $H$ is {\sf R}-bounded
as it is easy to check.
\end{remark}

\medskip
We have the following result which is, among others, an extension of
Hernandez's result \cite[Theorem 5.1]{He} stating that $C_p(X)$ is
strictly {\sf M}-bounded if and only if $X$ is pseudocompact.

\begin{theorem} \label{cp-m-bounded} Let $X$ be a Tychonoff space. The following are equivalent:
\begin{itemize}
\item[{\rm (1)}] $C_p(X)$ is strictly {\sf H}-bounded;

\item[{\rm (2)}] $C_p(X)$ is {\sf H}-bounded;

\item[{\rm (3)}] $C_p(X)$ is strictly {\sf M}-bounded;

\item[{\rm (4)}] $C_p(X)$ is {\sf M}-bounded;

\item[{\rm (5)}]  $X$ is pseudocompact;

\item[{\rm (6)}] $C_m(X)$ is first countable;

\item[{\rm (7)}] $C_m(X)$ is strictly Fr\'echet-Urysohn;

\item[{\rm (8)}] $C_m(X)$ has countable strong fan tightness;

\item[{\rm (9)}] $C_m(X)$ has countable fan tightness;

\item[{\rm (10)}] $C_m(X)$ has countable tightness; %%4 -9 equivalent by \cite{McG} m%

\item[{\rm (11)}] $C_m(X)$ is completely metrizable; %% 10 - 12 and 4 in \cite{hola-zsilinsky}

\item[{\rm (12)}] $C_m(X)$ is \v Cech-complete;

\item[{\rm (13)}] $C_m(X)$ is hereditarily Baire.
\end{itemize}
\end{theorem}
$\proof$ (1) $\Rightarrow$ (2) $\Rightarrow$ (4) and (1)
$\Rightarrow$ (3) $\Rightarrow$ (4) are trivial. (4) $\Rightarrow$
(5) was proved in \cite{He}.

\smallskip (5) $\Rightarrow$ (1) If $X$ is
pseudocompact , then $C_p(X)$ is a subgroup of a $\sigma$-compact
group $\bigcup_{n\in\mathbb N}[-n,n]^X$ and thus $C_p(X)$ is
strictly {\sf H}-bounded.

\smallskip
Other equivalences have been proved in \cite{DiMaio-Hola-Holy-McCoy,
McG, hola-zsilinsky}.
 $\eproof$

\begin{remark} \label{rem-to-th-cp-m-bounded} \rm Theorem
\ref{cp-m-bounded} remains true for the space $C_p(X,\mathbb R^n)$,
$n\in \mathbb N$. It follows from $C_p(X,\mathbb R^n) \cong
\prod_{i\le n}C_p(X)$ and the fact that (strict) $H$-boundedness is
finitely multiplicative. However, this theorem is not true for
$C_p(X,\mathbb R^{\omega})$. If $C_p(X,\mathbb R^{\omega})$ is {\sf
H}-bounded, then for any $x\in X$ its subspace $C_p(\{x\}, \mathbb
R^{\omega})$ must be also {\sf H}-bounded. But, $C_p(\{x\}, \mathbb
R^{\omega})$ is homeomorphic to $\mathbb R^{\omega}$ which is not
{\sf M}-bounded (see \cite{He}).
\end{remark}

\begin{theorem} \label{cm-cu-aleph-0-bounded} Let $X$ be a Tychonoff space. The following are
equivalent:
\begin{itemize}
\item[{\rm (1)}] $C_{\gamma}(X)$ is $\aleph_0$-bounded;

\item[{\rm (2)}] $C_m(X)$ is $\aleph_0$-bounded;

\item[{\rm (3)}] $C_u(X)$ is $\aleph_0$-bounded;

\item[{\rm (4)}] $C_u(X)$ is separable;

\item[{\rm (5)}] $X$ is compact and metrizable;

\item[{\rm (6)}] $C_{\gamma}(X)$ is second countable; %% 4-8 in \cite{hola-zsilinsky}

\item[{\rm (7)}] $C_{\gamma}(X)$ has a countable network;

\item[{\rm (8)}] $C_{\gamma}(X)$ is separable;

\item[{\rm (9)}] $C_{\gamma}(X)$ is ccc;

\item[{\rm (10)}] $C_m(X)$ is separable;

\item[{\rm (11)}] $C_m(X)$ is Lindel\"{o}ff;

\item[{\rm (12)}] $C_m(X)$ is ccc.

\end{itemize}
\end{theorem}
$\proof$ (1) $\Rightarrow$ (2) $\Rightarrow$ (3) is clear because if
$C_{\gamma}(X)$ is $\aleph_0$-bounded, then also $C_m(X)$ and
$C_u(X)$ are $\aleph_0$-bounded.

(3) $\Rightarrow$ (4) Since the group $C_u(X)$ is metrizable, it is
first countable. Let $\{U_n:n \in \mathbb N\}$ be a countable local
base at $f_0\in C_u(X)$. As $C_u(X)$ is $\aleph_0$-bounded, for each
$n \in \mathbb N$ there is a sequence $(f_{n,m}:m\in \mathbb N)$
such that $C_u(X) = \bigcup_{m\in\mathbb N}(f_{n,m} + U_n)$. It is
easy to check that the countable set $\{f_{n,m} : n,m \in \mathbb
N\}$ is dense in $C_u(X)$.

(4) $\Rightarrow$ (5) Since  $C_u(X)$ is separable and metrizable it
is second countable which implies that $X$ is compact and metrizable
\cite[Theorem 4.2.4]{MN}.

(5) $\Rightarrow$ (6) $\Rightarrow$ (7) $\Rightarrow$ (8)
$\Rightarrow$ (9) were proved in \cite{hola-zsilinsky}, while
equivalence of (5) with each of (10)--(12) can be found in
\cite{DiMaio-Hola-Holy-McCoy}.

(9) $\Rightarrow$ (1) follows from the fact that $C_{\gamma}(X)$ is
ccc and thus $\aleph_0$-bounded.
 $\eproof$

\medskip
{\bf Note.} Observe that $\tau_u = \tau_m= \tau_{\gamma}$ when $X$
is compact.

\medskip
The following two results describe boundedness properties of the
space $C_k(X)$. The first of these results may be compared with
\cite[Corollary 4.2.7]{MN} (which is a result from \cite{ntantu}),
but we give a new, direct proof of it.

\begin{theorem} \label{ck-aleph-0-bounded} {\rm (\cite{ntantu})} Let $X$ be a Tychonoff space. The following are equivalent:

\begin{itemize}
\item[{\rm (1)}] $C_k(X)$ is $\aleph_0$-bounded;

\item[{\rm (2)}]  every compact set in $X$ is metrizable.
\end{itemize}
\end{theorem}

$\proof$  (1) $\Rightarrow$ (2) Let $A$ be a compact set in $X$.
Then the mapping $\pi : C_k(X) \to C_k(A)$ defined by $\pi(f) =
f\restriction A$ is a continuous epimorphism. (Since $A$ is a
compact subset of a Tychonoff space $X$, every $f\in C(A)$ has a
continuous extension $f^{\ast} \in C(X)$.) Thus by Fact 2 $C_k(A) =
C_u(A)$ is $\aleph_0$-bounded. By Theorem 2.4 $A$ must be
metrizable.

\smallskip
(2) $\Rightarrow$ (1) Let $U$ be a $\tau_k$-open neighborhood of the
constantly zero function $f_0 : X \to \mathbb R$. There is a compact
set $A$ in $X$ and $\varepsilon > 0$ such that
\[
U \supset W(f_0, A, \varepsilon ) = \{f \in C(X) : |f(x) - f_0(x)| <
\varepsilon \mbox{ for each } x \in A\}.
\]
Since $A$ is compact and metrizable, $C_u(A) = C_k(A)$ is separable.
Let $D= \{f_i : i \in \mathbb N\}$ be a countable dense set in
$C_u(A)$. For each $i\in\mathbb N$ let $f_i^{\ast} : X \to \mathbb
R$ be a continuous extension of $f_i : A \to \mathbb R$. Put
$D^{\ast} = \{f_i^{\ast}:i \in \mathbb N\}$. We claim that
\[
C(X) \subset U + D^{\ast} = \{f + f_i^{\ast} : f \in U, i \in
\mathbb N\}.
\]
Let $g \in C(X)$. Since $\{f_i : i \in\mathbb N\}$ is dense in
$C_u(A)$ there is $i \in \mathbb N$ such that (setting $g_A =
g\restriction A$)
\[
|g_A(x) -  f_i(x)| < \varepsilon \mbox{ for each } x \in A.
\]
Thus $|g(x) - f_i^{\ast}(x)| < \varepsilon$ for each $x \in A$, i.e.
$g - f_i^{\ast} \in U$. Therefore,  $g = (g-f_i^{\ast})+f_i^{\ast}
\in U + D^{\ast}$.
 $\eproof$

%%Th 2.4
\begin{theorem} \label{c_k-m-bounded}  Let X be a Tychonoff space. The following are equivalent:
\begin{itemize}
\item[{\rm (1)}] $C_k(X)$ is strictly {\sf H}-bounded;

\item[{\rm (2)}]  $C_k(X)$ is {\sf H}-bounded;

\item[{\rm (3)}] $C_k(X)$ is strictly {\sf M}-bounded;

\item[{\rm (4)}] $C_k(X)$ is {\sf M}-bounded;

\item[{\rm (5)}]  $X$ is pseudocompact and every compact set in $X$ is finite.
\end{itemize}
\end{theorem}

$\proof$ (1) $\Rightarrow$ (2) $\Rightarrow$ (4) and (1)
$\Rightarrow$ (3) $\Rightarrow$ (4) are clear.

\smallskip
(4) $\Rightarrow$ (5) If $C_k(X)$ is {\sf M}-bounded, then also
$C_p(X)$ is {\sf M}-bounded. Thus by Theorem \ref{cp-m-bounded} $X$
is pseudocompact. Let $A$ be a compact set in $X$. The restriction
map $\pi : C_k(X) \to C_k(A)$ is a continuous epimorphism. By Fact
2, $C_k(A)$ is {\sf M}-bounded. Since $C_k(A)= C_u(A)$, by Theorem
\ref{cm-cu-m-bound} below, $A$ must be finite.

\smallskip
(5) $\Rightarrow$ (1) Since every compact set in $X$ is finite, we
have $C_k(X) = C_p(X)$. Thus by Theorem \ref{cp-m-bounded} we are
done.
 $\eproof$

\begin{remark} \label{rem-to-th-c_k-m-bounded} \rm By arguing as in
Remark \ref{rem-to-th-cp-m-bounded} we conclude that Theorem
\ref{c_k-m-bounded} is true for $C_k(X,\mathbb R^n)$, $n\in\mathbb
N$, but not for $C_k(X,\mathbb R^{\omega})$.
\end{remark}

\begin{theorem} \label{cm-cu-m-bound} Let $X$ be a Tychonoff space. The following are equivalent:
\begin{itemize}
\item[{\rm (1)}] $C_m(X)$ is $\sigma$-compact;

\item[{\rm (2)}] $C_m(X)$ is strictly {\sf H}-bounded;

\item[{\rm (3)}] $C_m(X)$ is {\sf H}-bounded;

\item[{\rm (4)}] $C_m(X)$ is strictly {\sf M}-bounded;

\item[{\rm (5)}] $C_m(X)$ is {\sf M}-bounded;

\item[{\rm (6)}] $C_u(X)$ is {\sf M}-bounded;

\item[{\rm (7)}] $X$ is finite;

\item[{\rm (8)}] $C_p(X)$ is $\sigma$-compact.
\end{itemize}
\end{theorem}

$\proof$  (6) $\Rightarrow$ (7) {\sf M}-boundedness of $C_u(X)$
implies $\aleph_0$-boundedness of $C_u(X)$. By Theorem
\ref{cm-cu-aleph-0-bounded} $X$ must be compact and metrizable. Thus
$C_u(X)$ is a Banach space. By \cite{Ba} an infinite-dimensional
Banach space (considered as an Abelian topological group) admits no
isomorphic embedding into a product of {\sf M}-bounded groups. Thus
$X$ must be finite.

\smallskip
(7) $\Leftrightarrow$  (8) is known \cite{Ar}.
 $\eproof$

\begin{remark} \rm The above theorem is true also for $C(X,\mathbb
R^n)$, $n\in \mathbb N$, but not for $C(X,\mathbb R^{\omega})$.
Observe also that we can extend the list in this theorem by the
corresponding properties of the space $C_{\gamma}(X)$.
\end{remark}

\medskip
The following theorem (and the note after it) is a variant of a
result from \cite{Bk}.

\begin{theorem} Let $X$ be a Tychonoff space and $G$ a subgroup of $C_u(X)$. The following are
equivalent:
\begin{itemize}
\item[{\rm  (1)}] $G$ is strictly {\sf R}-bounded subgroup of $C_u(X)$;

\item[{\rm (2)}] $G$ is countable.
\end{itemize}
\end{theorem}
$\proof$ Only (1) implies (2) need the proof. We use the fact that
$C_u(X)$ is metrizable by the left invariant sup-metric $d$, so that
$G$ is so. Let $(U_n : n \in \mathbb N)$ be a countable local base
of the identity $f_0 \in G$ such that ${\rm diam}_d(U_n) < 1/2^{n}$,
$n\in \mathbb N$. Let $\varphi$ be a strategy by TWO. For each $U_n$
TWO picks a point $h_n = \varphi(U_n) \in G$. Let $H_0 =
\bigcap_{n\in \mathbb N}(h_n+U_n)$.   For a given finite sequence
$n_1,\cdots, n_k$ in $\mathbb N$, define $H_{n_1,\cdots,n_k} =
\bigcap_{n\in\mathbb N}\varphi(U_{n_1}, \cdots, U_{n_k},U_n) + U_n$.

Let $S$ denote the set of all finite sequences in $\mathbb N$.

\smallskip
{\bf Claim 1.}  $G\subset \bigcup_{s \in S}H_{s}$.

Suppose, on the contrary, that there is $f\in G\setminus
\bigcup_{s\in S}H_s$. Then there is $n_1\in\mathbb N$ such that $f
\notin \varphi(U_{n_1}) + U_{n_1}$. Further, there is $n_2\in\mathbb
N$ such that $f\notin \varphi(U_{n_1}, U_{n_2}) + U_{n_2}$. And so
on. In this way we obtain a sequence $n_1, n_2, \cdots, n_k, \cdots$
in $\mathbb N$ and a sequence $U_{n_1}, \cdots, U_{n_k}, \cdots$ of
neighborhoods of $f_0$ such that $f\notin \varphi(U_{n_1},\cdots,
U_{n_k}) +U_{n_k}$. This means that there is a play according to the
strategy $\varphi$ lost by TWO, which is a contradiction.

\smallskip
{\bf Claim 2.} Each $H_{s}$ has at most one element.

Because the metric $d$ on $G$ is left-invariant, we have ${\rm
diam}_d(\varphi(U_{n_1}, \cdots, U_{n_k}, U_n) = {\rm diam}_d(U_n)$,
so that by assumption on the base $\{U_n:n\in\mathbb N\}$ of $f_0$,
we conclude that the diameters of $H_{s}$ tends to $0$. $\eproof$

\begin{note} \rm The above theorem remains true for a subgroup $H$ of a
topological group $G$ metrizable by a left-invariant metric.
\end{note}

%\begin{proposition} If $G$ is a {\sf GN}-bounded group metrizable by a left-invariant metric,  then
%$G$ is countable.
%\end{proposition}
%$\proof$ Let $(U_n:n\in\mathbb N)$ be a countable neighbourhood base
%for $e\in G$ such that ${\rm diam}_d(U_n) < 1/2^n$. Since $G$ is
%{\sf GN}-bounded for each $n\in \mathbb N$ there is $x_n\in G$ such
%that $G = \bigcup_{n\in \mathbb N}\bigcap_{m > n}(x_m + U_m)$. For
%each $n$ the set $\bigcap_{m > n}(x_m + U_m)$ contains one element,
%so that $G$ is countable. $\eproof$

%\begin{corollary} If $X$ is a Tychonoff space, then each {\sf
%GN}-bounded subgroup of $C_u(X)$ is countable.
%\end{corollary}

%Note that by Fact 5 the group $C_u(X)$ is never {\sf GN}-bounded.

\section{Examples}

\begin{example} \rm (1) The ordinal space $[0,\omega_1)$ is
 pseudocompact, and thus, by Theorem \ref{cp-m-bounded},
 $C_p([0,\omega_1))$ is $H$-bounded. On the other hand,
 $C_k([0,\omega_1))$ is not $H$-bounded by Theorem
 \ref{c_k-m-bounded}. Also,
 $C_m([0,\omega_1)) = C_u([0,\omega_1))$ is not $H$-bounded, even it is not
 $\aleph_0$-bounded by Theorem \ref{cm-cu-aleph-0-bounded}. Note that the
 compact-open topology on $C([0,\omega_1))$ is strictly weaker
 than the topology of uniform convergence.

\smallskip
(2) $C_k(\mathbb R)$ and $C_k([0,1])$ are $\aleph_0$-bounded
(Theorem \ref{ck-aleph-0-bounded}), but not {\sf M}-bounded (Theorem
\ref{c_k-m-bounded}). The spaces $C_u(\mathbb R)$ and $C_m(\mathbb
R)$ also are not {\sf M}-bounded (Theorem \ref{cm-cu-m-bound}). The
spaces $C_u([0,1]) = C_m([0,1])= C_{\gamma}([0,1])$ are not {\sf
M}-bounded.

\smallskip
(3) Let $T_{\infty}$ be the deleted Tychonoff plank $[0,\omega_1]
\times [0,\omega] \setminus \{(\omega_1,\omega)\}$ (see
\cite{steen-seebach}).  This space is pseudocompact but not compact,
so that $C_p(T_{\infty})$ is $\aleph_0$-bounded, but neither
$C_k(T_{\infty})$ nor $C_u(T_{\infty}) = C_m(T_{\infty})$ is
$\aleph_0$-bounded.

\smallskip
(4) For each infinite subset $A$ of $\mathbb N$ pick a point $x_A
\in {\rm Cl}_{\beta \mathbb N}(A) \setminus \mathbb N$ and set $X =
\mathbb N \cup \{x_A: A \mbox{ is an infinite subset of } \mathbb
N\}$. This space is pseudocompact and each compact subset of $X$ is
finite (see  \cite[Example 2.6]{kundu-garg}). By Theorem
\ref{cp-m-bounded} $C_k(X) = C_p(X)$ is {\sf H}-bounded. But,
$C_m(X)$ is not {\sf H}- bounded.

\smallskip
(5) A family $\mathcal A$ of infinite subsets of $\mathbb N$ is
almost disjoint if for all distinct $A,B \in \mathcal A$ the set
$A\cap B$ is finite. An almost disjoint family is maximal if it is
not contained in any other almost disjoint family. A topological
space $X$ is a Mr\'owka-Isbell space or a $\Psi$-space if it has the
form $\Psi(\mathcal A) = \mathbb N \cup \mathcal A$, where $\mathcal
A$ is an almost disjoint family, and its topology is generated by
the following base: each $\{n\} \subset \mathbb N$ is open, and an
open canonical neighborhood of $A \in \mathcal A$ is of the form
$\{A\} \cup (A\setminus F)$, where $F \subset \mathbb N$ is finite
\cite{engelking}. It is known that $\Psi(\mathcal A)$ is
pseudocompact if and only if $\mathcal A$ is maximal almost disjoint
family. So, if $\mathcal A$ is a maximal almost disjoint family,
then by Theorem \ref{cp-m-bounded} $C_p(\Psi(\mathcal A))$ is {\sf
M}-bounded, but, by Theorems \ref{c_k-m-bounded} and
\ref{cm-cu-m-bound}, neither of $C_k(\Psi(\mathcal A))$,
$C_u(\Psi(\mathcal A))$, $C_m(\Psi(\mathcal A))$ and
$C_{\gamma}(\Psi(\mathcal A))$ is {\sf M}-bounded.
\end{example}

\section*{Acknowledgements}

The authors are grateful to referees for several comments and
suggestions which led to the improvement of the exposition. The
first author would like to thank the support of Vega 2/0006/16.

\footnotesize{
 }

\end{document}